\documentclass[11pt]{article}

\usepackage[margin=1in]{geometry}
\usepackage[T1]{fontenc}
\usepackage{lmodern}
\usepackage{amsmath,amssymb,amsfonts,amsthm,mathtools,mathrsfs}
\usepackage{aliascnt}
\usepackage{enumitem}
\usepackage{xcolor}
\usepackage{hyperref}
\usepackage{microtype}
\usepackage{authblk}
\usepackage[nameinlink,capitalise]{cleveref}

\hypersetup{colorlinks=true,linkcolor=blue!50!black,citecolor=blue!50!black,urlcolor=blue!50!black}

\numberwithin{equation}{section}

\theoremstyle{plain}
\newtheorem{theorem}{Theorem}[section]
\newaliascnt{lemma}{theorem}
\newtheorem{lemma}[lemma]{Lemma}
\aliascntresetthe{lemma}
\newaliascnt{proposition}{theorem}
\newtheorem{proposition}[proposition]{Proposition}
\aliascntresetthe{proposition}
\newaliascnt{corollary}{theorem}
\newtheorem{corollary}[corollary]{Corollary}
\aliascntresetthe{corollary}
\newaliascnt{conjecture}{theorem}
\newtheorem{conjecture}[conjecture]{Conjecture}
\aliascntresetthe{conjecture}
\newaliascnt{problem}{theorem}

\aliascntresetthe{problem}

\theoremstyle{definition}
\newaliascnt{definition}{theorem}
\newtheorem{definition}[definition]{Definition}
\aliascntresetthe{definition}

\theoremstyle{remark}
\newaliascnt{remark}{theorem}
\newtheorem{remark}[remark]{Remark}
\aliascntresetthe{remark}

\crefname{theorem}{Theorem}{Theorems}
\Crefname{theorem}{Theorem}{Theorems}
\crefname{lemma}{Lemma}{Lemmas}
\Crefname{lemma}{Lemma}{Lemmas}
\crefname{proposition}{Proposition}{Propositions}
\Crefname{proposition}{Proposition}{Propositions}
\crefname{corollary}{Corollary}{Corollaries}
\Crefname{corollary}{Corollary}{Corollaries}
\crefname{definition}{Definition}{Definitions}
\Crefname{definition}{Definition}{Definitions}
\crefname{remark}{Remark}{Remarks}
\Crefname{remark}{Remark}{Remarks}
\crefname{conjecture}{Conjecture}{Conjectures}
\Crefname{conjecture}{Conjecture}{Conjectures}
\crefname{problem}{Problem}{Problems}
\Crefname{problem}{Problem}{Problems}
\crefname{section}{Section}{Sections}
\Crefname{section}{Section}{Sections}
\crefname{equation}{Equation}{Equations}
\Crefname{equation}{Equation}{Equations}

\DeclareMathOperator{\ex}{ex}
\newcommand{\F}{\mathbb{F}}
\newcommand{\A}{\mathbb{A}}
\newcommand{\bP}{\mathbb{P}}

\begin{document}
\baselineskip=16pt
\title{Rational Tur\'an exponents of a single $k$-uniform hypergraph\footnote{Corresponding author: Chang Liu (liuchang\_@nudt.edu.cn)}}
\author{Chang Liu, Jue Dai}
\affil{College of Science, National University of Defense Technology,\linebreak Changsha 410073, China}
\date{}
\maketitle

\begin{abstract}
	For a fixed $k$-uniform hypergraph $F$, let $\ex_k(n,F)$ denote the maximum number of edges in an $F$-free $k$-uniform hypergraph on $n$ vertices. We prove that for every integer $k\ge2$ and every rational $r\in[0,1]$ there exists a $k$-uniform hypergraph $F_{k,r}$ such that
	\begin{equation*}
		\ex_k(n,F_{k,r})=\Theta(n^{k-r}).
	\end{equation*}
	The case $k=2$ is the rational exponents conjecture of Erd\H{o}s and Simonovits \cite{Erdos1988}, also known as Erd\H{o}s problem \#571 \cite{BloomErdoswebsite}; a proof was announced in September 2026 and is expounded by Conlon et al.\ \cite{Conlon-Janzer-Jiang-Kim-Lee-Liu-Longbrake-Yepremyan}. For $k\ge3$ the theorem extends their result to $k$-uniform hypergraphs. In the single-hypergraph setting, this strengthens Fitch's finite-family theorem \cite{Fitch} by realizing each exponent with a single forbidden $k$-uniform hypergraph for $0\le r\le1$. We also show, via an explicit two-edge construction, that for $k\ge4$ this strengthening fails for $\lceil k/2\rceil<r\le k-1$: no single $k$-uniform hypergraph $F$ then satisfies $\ex_k(n,F)=\Theta(n^{k-r})$. The main result of this paper has been formalized in Lean 4, conditional on two explicitly stated axioms.
\end{abstract}



\section{Introduction}

A $k$-uniform hypergraph, or $k$-graph, is a pair $H=(V(H),E(H))$ in which $E(H)$ is a set of $k$-element subsets of $V(H)$. The degree $d_H(v)$ of a vertex $v$ is the number of edges containing $v$. For $k$-graphs $F$ and $H$, a copy of $F$ in $H$ is an injective map $f\colon V(F)\to V(H)$ such that $f(e)\in E(H)$ for every $e\in E(F)$. A $k$-graph $H$ is $F$-free if it contains no copy of $F$, and
\begin{equation*}
	\ex_k(n,F)=\max\{e(H):|V(H)|=n,\ H\text{ is }F\text{-free}\}.
\end{equation*}
When $k=2$, we write $\ex(n,F)$ for $\ex_2(n,F)$. We consider only finite families $\mathscr{F}$ of $k$-graphs. A $k$-graph $H$ is $\mathscr{F}$-free if it is $F$-free for every $F\in\mathscr{F}$, and $\ex_k(n,\mathscr{F})$ is defined by the same maximum taken over $\mathscr{F}$-free $k$-graphs. A $k$-graph is $k$-partite if its vertex set is the disjoint union of $k$ classes $V_1,\dots,V_k$ and every edge contains exactly one vertex of each class.

The Erd\H{o}s--Stone--Simonovits theorem \cite{ErdosSimonovits1966,ErdosStone1946} gives
\begin{equation*}
	\ex(n,F)=\left(1-\frac{1}{\chi(F)-1}+o(1)\right)\binom n2
\end{equation*}
for every fixed graph $F$ with $\chi(F)\ge3$. For bipartite $F$, this theorem gives only $\ex(n,F)=o(n^2)$, whereas the K\H{o}v\'ari--S\'os--Tur\'an theorem gives $\ex(n,F)=O(n^{2-\delta})$ for some $\delta>0$ \cite{KST}. Determining the sharp exponent remains difficult even for $K_{4,4}$ and $C_8$. These difficulties led to general questions about the possible orders of magnitude of $\ex(n,F)$. In particular, Erd\H{o}s conjectured that every bipartite graph has a rational exponent
\cite{ErdosSimonovits1969}, and Erd\H{o}s and Simonovits conjectured, conversely, that every rational $\alpha\in[1,2]$ is realized by a finite
graph family $\mathscr{F}$ with $\ex(n,\mathscr{F})=(c_{\mathscr{F}}+o(1))n^\alpha$ for some $c_{\mathscr{F}}>0$ \cite{Erdos1981,Erdos1988}.

Bukh and Conlon \cite{BukhConlon2018} proved the order-of-magnitude form of this family conjecture: every rational $\alpha\in[1,2]$ satisfies $\ex(n,\mathscr{F})=\Theta(n^\alpha)$ for some finite graph family $\mathscr{F}$. Their proof uses the random algebraic method introduced by Bukh \cite{Bukh2015}. Fitch extended this result to uniform hypergraphs \cite{Fitch}.

\begin{theorem}[Fitch \cite{Fitch}]
	For every integer $k\ge2$ and every rational $r\in[0,k-1]$, there exists a finite family $\mathscr F$ of $k$-uniform hypergraphs such that
	\begin{equation*}
		\ex_k(n,\mathscr F)=\Theta(n^{k-r}).
	\end{equation*}
\end{theorem}

The corresponding assertion for a single forbidden graph is the rational exponents conjecture of Erd\H{o}s and Simonovits: for every rational $\alpha\in[1,2]$, there is a graph $F$ with $\ex(n,F)=\Theta(n^\alpha)$ \cite{Erdos1988}. The endpoint $\alpha=2$ is realized by a non-bipartite graph; the substantive assertion concerns bipartite graphs and $1\le\alpha<2$. The distinction between one graph and a finite family is essential. The Bukh--Conlon construction \cite{BukhConlon2018} supplies lower bounds for large powers of balanced rooted graphs, but the corresponding elementary counting argument forces only a member of a finite family, since different copies may overlap outside the roots.

Work on single forbidden graphs addressed this intersection problem through increasingly refined embedding arguments. Kang, Kim and Liu \cite{KangKimLiu} realized $2-a/b$ whenever $0<a<b$ and $b\equiv\pm1\pmod a$, and showed that their subdivision conjecture would imply the full rational exponents conjecture. Further progress included the subdivision results of Conlon, Janzer and Lee \cite{ConlonJanzerLee2021}, the realization of $1+a/b$ for $b>a^2$ by Jiang and Qiu \cite{JiangQiu2023}, and that of $2-a/b$ for $b>a$ and $b\ge(a-1)^2$ by Conlon and Janzer \cite{ConlonJanzer2022}. In September 2026, a proof produced with the model GPT-6 Astra was announced in connection with Epoch AI's benchmarking (see \cite{BloomErdoswebsite,Conlon-Janzer-Jiang-Kim-Lee-Liu-Longbrake-Yepremyan}). Conlon, Janzer, Jiang, Kim, Lee, Liu, Longbrake and Yepremyan subsequently presented the argument in \emph{A proof of the rational exponents conjecture} \cite{Conlon-Janzer-Jiang-Kim-Lee-Liu-Longbrake-Yepremyan}, explaining its relationship to the preceding work. Their exposition presents a complete proof of the single-graph conjecture.

It is therefore natural to ask whether Fitch's finite family can likewise be replaced by a single $k$-uniform hypergraph: for which rational $r\in[0,k-1]$ does there exist $F$ with $\ex_k(n,F)=\Theta(n^{k-r})$?

The main result of this paper gives such a realization for $0\le r\le1$.

\begin{theorem}\label{thm:mainresult-1}
For every integer $k\ge2$ and every rational $r\in[0,1]$, there exists a $k$-uniform hypergraph $F_{k,r}$ such that
\begin{equation*}
	\ex_k(n,F_{k,r})=\Theta(n^{k-r}).
\end{equation*}
\end{theorem}

For $r=0$, we take $F_{k,0}=K_{k+1}^{(k)}$, the complete $k$-uniform hypergraph on $k+1$ vertices, and prove the claim in \cref{sec:proof}; for $r>0$ the hypergraph $F_{k,r}$ is constructed in \cref{sec:proof} as a rooted power of a suspension.

The upper bound in \cref{thm:mainresult-1} is an immediate consequence of the suspension inequality of Cheng, Gerbner, Hama Karim and Zhou \cite{cheng2026turan}. The lower bound is more delicate. It requires a dense $k$-partite $k$-graph $X$ avoiding one fixed rooted power of the suspended graph; equivalently, for every $(k-2)$-set $S\subseteq V(X)$, the link $\{e\setminus S:e\in E(X),\ S\subseteq e\}$ must avoid the corresponding rooted power of the base graph. We obtain such an $X$ by a random algebraic construction over a finite field.

The proof uses two main ingredients. The first is a uniform dichotomy for the number of $\F_q$-points of constructible sets of bounded complexity. It combines Heintz's elimination bound \cite{Heintz}, B\'ezout's inequality \cite{Fulton}, and the Lang--Weil estimate in the explicit form of Cafure--Matera \cite{CafureMatera}. The second is a moment estimate for the number of copies of the suspended graph above a fixed root set, which follows from the balance condition.

We also record a limitation of the single-hypergraph approach.

\begin{proposition}\label{prop:obstruction}
Let $k\ge4$ and let $r$ be a rational number with $\lceil k/2\rceil<r\le k-1$. Then there is no $k$-uniform hypergraph $F$ with $\ex_k(n,F)=\Theta(n^{k-r})$.
\end{proposition}

The proof is given in Section~\ref{sec:obstruction}. For $k=2$, \cref{thm:mainresult-1} settles the single-graph case completely. For $k\ge3$, our results leave the interval $1<r\le\lceil k/2\rceil$ unresolved. We discuss the subinterval $1<r\le k/2$ in Section~\ref{sec:concluding-remark}.


\section{Preliminaries}\label{sec:prelim}

\subsection{Rooted graphs and the ordinary-graph input}

A rooted graph is a pair $(T,R)$ in which $T$ is a graph and $\emptyset\ne R\subsetneq V(T)$. The vertices of $R$ are the roots. We write $T$ for $(T,R)$ when the root set is clear. For $S\subseteq V(T)\setminus R$, let $e_T(S)$ be the number of edges of $T$ incident with at least one vertex of $S$, and put
\begin{equation*}
	\rho_T(S)=\frac{e_T(S)}{|S|},\qquad \rho(T)=\rho_T(V(T)\setminus R).
\end{equation*}
We say that $T$ is balanced if $\rho_T(S)\ge\rho(T)$ for every non-empty $S\subseteq V(T)\setminus R$. We write $\rho(T)=b/a$, where $a,b$ are positive integers with $a\le b$.

For $p\ge1$, the $p$-th power $T^{(p)}$ is obtained from $p$ vertex-disjoint copies of $T$ by identifying the $p$ copies of each root. Thus the copies are disjoint outside $R$.

We use the following rooted form of the Erd\H{o}s--Simonovits theorem.

\begin{proposition}[{\cite{BukhConlon2018,Conlon-Janzer-Jiang-Kim-Lee-Liu-Longbrake-Yepremyan}}]\label{prop:rooted}
For every rational $r=a/b\in(0,1]$, there exists a finite connected bipartite rooted graph $(T,R)$ with $V(T)\setminus R\ne\emptyset$ and $\rho(T)=b/a$ such that
\begin{enumerate}
\item[\textnormal{(i)}] $b|S|\le a e_T(S)$ for every non-empty $S\subseteq V(T)\setminus R$;
\item[\textnormal{(ii)}] $\ex(n,T^{(p)})=O(n^{2-a/b})$ for every fixed $p\ge1$.
\end{enumerate}
\end{proposition}

For $r=1$, the rooted edge $K_2$ satisfies both assertions. In general, \cite[Proposition~2.3]{Conlon-Janzer-Jiang-Kim-Lee-Liu-Longbrake-Yepremyan} gives a balanced connected bipartite rooted graph $T$ with $\rho(T)=b/a$, and equation~(4) in Section~2.4 of that paper gives (ii).

\subsection{Suspension hypergraphs}

Let $F$ be a graph and $k\ge2$. The $k$-uniform suspension $\mathcal S^kF$ is the $k$-graph obtained from $F$ by adding a set $W$ of $k-2$ new vertices to every edge. Thus
\begin{equation*}
	V(\mathcal S^kF)=V(F)\sqcup W,\qquad E(\mathcal S^kF)=\{e\cup W:e\in E(F)\}.
\end{equation*}
When $k=2$, $W=\emptyset$ and $\mathcal S^2F=F$. For a rooted graph $T$, the suspension $\mathcal S^kT$ is rooted on $R\cup W$.

The following inequality is due to Cheng, Gerbner, Hama Karim and Zhou \cite{cheng2026turan}.

\begin{lemma}[{\cite{cheng2026turan}}]\label{lem:suspension-upper}
Let $k\ge2$ and let $F$ be a graph. Then for all $n$,
\[
\ex_k(n,\mathcal S^kF)
\le\binom{n}{k-2}\frac{\ex(n-k+2,F)}{\binom{k}{2}}
=(1+o(1))\frac{2\binom nk\ex(n,F)}{n^2}.
\]
\end{lemma}

Let $r=a/b\in(0,1]$ and let $T$ be as in \cref{prop:rooted}. For every fixed $p\ge1$, \cref{lem:suspension-upper} and \cref{prop:rooted}(ii) give
\begin{equation*}
\ex_k\bigl(n,\mathcal S^k(T^{(p)})\bigr)
\le\binom{n}{k-2}\frac{\ex(n-k+2,T^{(p)})}{\binom{k}{2}}
=O\bigl(n^{k-a/b}\bigr).
\end{equation*}
This is the upper bound in \cref{thm:mainresult-1} for $r>0$. The lower bound is proved in Section~\ref{sec:proof}.

\subsection{Algebraic geometry over finite fields}\label{sec:ag}

The lower bound in Section~\ref{sec:proof} uses a random algebraic construction. We collect the three effective results on which it depends: Heintz's elimination bound \cite{Heintz}, the Lang--Weil estimate \cite{LangWeil}, and the explicit estimate of Cafure--Matera \cite{CafureMatera}. We then derive a uniform dichotomy for bounded-complexity constructible sets. We use standard notions of dimension, degree and irreducibility as in \cite{Shafarevich}, and B\'ezout's inequality as in \cite{Fulton}. For an algebraic set, its degree is the sum of the degrees of all its geometric irreducible components, including components of every dimension. All degree bounds below use this convention.

\subsubsection{Complexity and Frobenius}

Let $q$ be a prime power and let $\overline{\F}_q$ be an algebraic closure of $\F_q$. Write $\A^N=\overline{\F}_q^N$ and $\bP^N$ for affine and projective $N$-space. An algebraic set $V\subseteq\A^N$ is defined over $\F_q$ if it is the zero set of polynomials with coefficients in $\F_q$. The $q$-power Frobenius map is
\begin{equation*}
	\Phi_q\colon\A^N\to\A^N,\qquad (x_1,\dots,x_N)\mapsto(x_1^q,\dots,x_N^q).
\end{equation*}
A set $V$ is defined over $\F_q$ if and only if $\Phi_q(V)=V$; in that case $V(\F_q)=\{x\in V:\Phi_q(x)=x\}$. This is standard Galois descent for Zariski-closed sets.

\begin{definition}\label{def:complexity}
A set $Y\subseteq\A^N$ has complexity at most $(N,m,D)$ if there are polynomials $g_1,\dots,g_u,h_1,\dots,h_v\in\F_q[x_1,\dots,x_N]$ with $u+v\le m$ and $\deg g_i,\deg h_j\le D$ such that
\[
Y=\{x\in\A^N:g_1(x)=\cdots=g_u(x)=0,\ h_1(x)\cdots h_v(x)\ne0\}.
\]
\end{definition}

\subsubsection{Effective B\'ezout bounds}

\begin{theorem}[Heintz {\cite{Heintz}}]\label{thm:heintz}
For all $N,m,D$ there exist integers $B_0,\Gamma_0\ge1$, depending only on $(N,m,D)$, such that every algebraic set $V\subseteq\A^N$ defined by at most $m$ polynomials of degree at most $D$ has at most $B_0$ irreducible components, and the sum of their degrees is at most $\Gamma_0$. The bounds are uniform in the coefficients.
\end{theorem}

\begin{proposition}[B\'ezout's inequality \cite{Fulton}]\label{prop:bezout}
Let $X,X'\subseteq\bP^N$ be irreducible projective varieties. Then
\begin{equation*}
	\deg(X\cap X')\le\deg X\cdot\deg X',
\end{equation*}
where the degree of the intersection is the sum of the degrees of its irreducible components, and $X\cap X'$ has at most $\deg X\cdot\deg X'$ irreducible components. The same bounds hold for affine varieties.
\end{proposition}

The proposition is the projective B\'ezout inequality; see \cite[Section~12.3]{Fulton}. For affine varieties, intersecting with the hyperplane at infinity can only add components.

\begin{corollary}\label{cor:descent}
Let $B_0,\Gamma_0$ be as in \cref{thm:heintz}, and define $\Gamma_{t+1}=\Gamma_t^2$ and $B_{t+1}=B_t\Gamma_t^2$ for $t\ge0$. Let $V\subseteq\A^N$ be defined over $\F_q$ and have at most $B_0$ irreducible components of degree at most $\Gamma_0$. For a component $C$ of $V$, form the tree whose root is $C$ and whose children of a non-Frobenius-stable node $G$ are the irreducible components of $G\cap\Phi_q(G)$. Then at depth $t$ the forest has at most $B_t$ nodes, each of degree at most $\Gamma_t$, and every branch has length at most $N$.
\end{corollary}

\begin{proof}
Induct on $t$. A non-stable node $G$ satisfies $G\ne\Phi_q(G)$. Since $\Phi_q$ is the $q$-power automorphism of $\overline{\F}_q$, it preserves irreducibility, dimension and degree; in particular, $\Phi_q(G)$ is irreducible, $\dim\Phi_q(G)=\dim G$ and $\deg(\Phi_q(G))=\deg G$. Thus every component of $G\cap\Phi_q(G)$ has dimension at most $\dim G-1$. By \cref{prop:bezout}, there are at most $(\deg G)^2\le\Gamma_t^2$ such components, each of degree at most $\Gamma_t^2=\Gamma_{t+1}$. The assertions follow. Since each non-stable step strictly lowers the dimension, every branch has length at most $N$; hence the tree is finite and all bounds depend only on $N,m,D$.
\end{proof}

\subsubsection{Point counts on varieties}

\begin{theorem}[Lang--Weil {\cite{LangWeil}}]\label{thm:lang-weil}
Let $V\subseteq\A^N$ be a geometrically irreducible variety defined over $\F_q$, of dimension $d\ge1$ and degree $\delta$. Then
\[
\bigl|\,|V(\F_q)|-q^d\,\bigr|
\le(\delta-1)(\delta-2)q^{d-1/2}+Cq^{d-1},
\]
where $C=C(N,\delta)$ depends only on $N$ and $\delta$.
\end{theorem}

\begin{theorem}[Cafure--Matera {\cite{CafureMatera}}]\label{thm:cafure-matera}
Let $V$ be as in \cref{thm:lang-weil}. If $q>2(d+1)\delta^2$, then
\begin{equation*}
	\bigl|\,|V(\F_q)|-q^d\,\bigr|\le(\delta-1)(\delta-2)q^{d-1/2}+5\delta^{13/3}q^{d-1}.
\end{equation*}
\end{theorem}

\begin{corollary}\label{cor:degree-uniform}
For all $N,\Gamma$ there exist $A=A(N,\Gamma)$ and $q_1=q_1(N,\Gamma)$ such that every geometrically irreducible variety $V\subseteq\A^N$ defined over $\F_q$, of dimension $d\ge1$ and degree at most $\Gamma$, satisfies
\begin{equation*}
	\bigl|\,|V(\F_q)|-q^d\,\bigr|\le Aq^{d-1/2}
\end{equation*}
for every $q\ge q_1$. Explicitly, $A=6\Gamma^{13/3}$ and $q_1=2(N+1)\Gamma^2+1$ suffice.
\end{corollary}

\begin{proof}
Since $d\le N$ and $\delta\le\Gamma$, the hypothesis $q\ge q_1$ implies $q>2(d+1)\delta^2$. The claim follows from \cref{thm:cafure-matera} and $q^{d-1}\le q^{d-1/2}$.
\end{proof}

\begin{lemma}[Affine degree bound \cite{CafureMatera}]\label{lem:affine-degree}
Let $Z\subseteq\A^N$ be an algebraic set of dimension $e$ and degree $\delta$. Then
\[
|Z(\F_q)|\le\delta q^e
\]
for every $q$.
\end{lemma}

This is the elementary affine degree bound; see \cite[Lemma~2.1]{CafureMatera}.

\begin{corollary}[Open Lang--Weil]\label{cor:open-lang-weil}
Let $V\subseteq\A^N$ be geometrically irreducible and defined over $\F_q$, of dimension $d\ge1$ and degree at most $\Gamma$. Let $S\subseteq V$ be a proper closed subset with $\deg(V\cap S)\le\Gamma'$. Then, for every $q\ge q_1(N,\Gamma)$,
\begin{equation*}
	|(V\setminus S)(\F_q)|\ge q^d-Aq^{d-1/2}-\Gamma'q^{d-1},
\end{equation*}
where $A=A(N,\Gamma)$ is as in \cref{cor:degree-uniform}. In particular, one may take
\[
q_0=\max\left\{q_1(N,\Gamma),\left\lceil(2A+2\Gamma')^2\right\rceil\right\},
\]
so that $|(V\setminus S)(\F_q)|\ge q/2$ for all $q\ge q_0$.
\end{corollary}

\begin{proof}
The first inequality follows from \cref{thm:lang-weil} and \cref{lem:affine-degree}, applied to $V\cap S$, whose dimension is at most $d-1$. The displayed choice of $q_0$ implies $Aq^{-1/2}+\Gamma'q^{-1}\le1/2$, so the right-hand side is at least $q^d/2\ge q/2$.
\end{proof}

\subsubsection{A uniform dichotomy for constructible sets}

\begin{theorem}\label{thm:uniform-dichotomy}
For all $N,m,D$ there exist $P=P(N,m,D)$ and $q_0=q_0(N,m,D)$ such that, for every $q\ge q_0$ and every $Y\subseteq\A^N$ of complexity at most $(N,m,D)$,
\begin{equation*}
	|Y(\F_q)|\le P\qquad\text{or}\qquad |Y(\F_q)|\ge q/2.
\end{equation*}
The constants are uniform in the coefficients of the defining polynomials.
\end{theorem}

\begin{proof}
Write $Y=V\setminus Z(h)$, where $V=Z(g_1,\dots,g_u)$, $h=h_1\cdots h_v$, and $u+v\le m$. Decompose $V$ over $\overline{\F}_q$ into geometrically irreducible components $C_1,\dots,C_s$. By \cref{thm:heintz}, $s\le B_0$ and $\deg C_\nu\le\Gamma_0$. Since $\deg h\le mD$, we have $\deg(C_\nu\cap Z(h))\le\deg C_\nu\deg h\le\Gamma_0\max\{1,mD\}$ for all $\nu$.

Discard the components contained in $Z(h)$, and for the remaining ones put $E_\nu=C_\nu\cap Z(h)$ and $W_\nu=C_\nu\setminus E_\nu$. Then $Y(\F_q)=\bigcup_\nu W_\nu(\F_q)$.

We claim that each $W_\nu(\F_q)$ is either small or large. If $C_\nu$ is Frobenius-stable and zero-dimensional, then $|W_\nu(\F_q)|\le\deg C_\nu\le\Gamma_0$. If $C_\nu$ is Frobenius-stable and positive-dimensional, then either $C_\nu\subseteq Z(h)$, in which case $W_\nu=\emptyset$, or \cref{cor:open-lang-weil} gives $|W_\nu(\F_q)|\ge q/2$.

If $C_\nu$ is not Frobenius-stable, form the descent tree of \cref{cor:descent}. At depth $t$ there are at most $B_t$ nodes of degree at most $\Gamma_t$ (the bound of \cref{cor:descent} applies to the forest, hence to each of its trees), and the depth is at most $N$. A node $G$ at depth $t$ not contained in $Z(h)$ satisfies
\[
\deg(G\cap Z(h))\le \deg G\,\deg h\le \Gamma_t\max\{1,mD\}.
\]
Thus all applications of \cref{cor:open-lang-weil} below use the uniform bounds $\Gamma_N$ and $\Gamma_N\max\{1,mD\}$. Call a node \emph{terminal} if it has no children. Every $x\in W_\nu(\F_q)$ lies in a terminal node: if $x$ lies in a non-stable node $G$, then $\Phi_q(x)=x\in G\cap\Phi_q(G)$, so $x$ lies in a child. In particular, a non-stable terminal node contains no $\F_q$-point. If a terminal node $G$ is positive-dimensional, Frobenius-stable, and not contained in $Z(h)$, then \cref{cor:open-lang-weil} gives $|G(\F_q)\setminus Z(h)|\ge q/2$, and hence $|W_\nu(\F_q)|\ge q/2$. Otherwise every terminal node is zero-dimensional, or is contained in $Z(h)$, or is non-stable. A non-stable terminal node contributes no $\F_q$-point; a zero-dimensional node contributes at most $\Gamma_N$ points; a node contained in $Z(h)$ contributes none. Extending terminal nodes to depth $N$ by unary nodes, we see that there are at most $B_N$ terminal nodes. Hence $|W_\nu(\F_q)|\le B_N\Gamma_N$.

Thus each $W_\nu(\F_q)$ is either at most $B_N\Gamma_N$ or at least $q/2$. If one is large, then $|Y(\F_q)|\ge q/2$. Otherwise $|Y(\F_q)|\le sB_N\Gamma_N\le B_0B_N\Gamma_N$. This proves the theorem with $P=B_0B_N\Gamma_N$ and $q_0=q_0(N,\Gamma_N,\Gamma_N\max\{1,mD\})$ from \cref{cor:open-lang-weil}.
\end{proof}

\begin{remark}\label{rem:dichotomy-uniform}
The constants $B_t$ and $\Gamma_t$ are obtained from $B_0,\Gamma_0$ by the fixed recursions of \cref{cor:descent}, and $B_0,\Gamma_0$ depend only on $N,m,D$ by \cref{thm:heintz}. Consequently $P$ and $q_0$ in \cref{thm:uniform-dichotomy} depend only on $(N,m,D)$, not on the coefficients of the defining polynomials. This uniformity is what allows the root-fibre constants $P_*,q_*$ in \cref{prop:root-fibre} to be independent of the host coefficients and of the root position.
\end{remark}

\subsubsection{Root fibres}

Fix a rooted $k$-partite $k$-graph $U$ with root set $R_U$, interior $A_U$, $h=|A_U|$ and $\varrho=|R_U|$. Fix $a,b\ge1$ and $D\ge1$. Let $X$ be a $k$-partite $k$-graph with vertex classes $V_1,\dots,V_k$, each identified with $\F_q^b$, whose edges are defined by polynomials $f_1,\dots,f_a$ of degree at most $D$: a $k$-tuple $(z_1,\dots,z_k)\in V_1\times\cdots\times V_k$ is an edge if and only if $f_1(z)=\cdots=f_a(z)=0$. Write $\kappa_U\colon V(U)\to[k]$ and $\kappa_X\colon V(X)\to[k]$ for the part maps.

A map $w\colon R_U\to V(X)$ is a \emph{legal root position} if it is injective and there is a permutation $\pi\in S_k$ such that
\begin{equation}\label{eq:legal-root}
	\kappa_X(w(r))=\pi(\kappa_U(r))\qquad\text{for every }r\in R_U.
\end{equation}
For a legal root position $w$, let $\Pi(w)$ be the set of permutations satisfying \eqref{eq:legal-root}, and let $C_w$ be the number of injective maps $\varphi\colon V(U)\to V(X)$ such that $\varphi|_{R_U}=w$, $\varphi(e)\in E(X)$ for every $e\in E(U)$, and $\kappa_X\circ\varphi=\pi\circ\kappa_U$ for some $\pi\in\Pi(w)$.

In the application $U$ has at least one edge. Then a map $\varphi$ determines $\pi$ on the $k$ vertices of any edge, hence on all of $[k]$; consequently $C_w$ counts embeddings, not pairs $(\varphi,\pi)$.

Equivalently, writing $C_w(\pi)$ for the number of embeddings compatible with a fixed $\pi$, we have
\[
C_w=\sum_{\pi\in S_k}C_w(\pi),
\]
where $C_w(\pi)=0$ for $\pi\notin\Pi(w)$: indeed, a compatible embedding has to satisfy \eqref{eq:legal-root} on every root. This is the form used in the counting argument below.

Thus a legal root position records one permutation of the parts, applied simultaneously to all roots; the embeddings counted by $C_w$ use a permutation from $\Pi(w)$. For the suspension $U$ constructed in Section~\ref{sec:proof}, the set $\Pi(w)$ is in fact a singleton: the $k-2$ vertices of $W$ fix $\pi$ on $k-2$ parts, and any root $r\in R$ fixes one of the two remaining parts. In the general discussion below we keep the notation $\Pi(w)$ to avoid treating this as a special case.

\begin{proposition}\label{prop:root-fibre}
There exist integers $P_*=P_*(U,a,b,D)$ and $q_*=q_*(U,a,b,D)$ such that, for every $q\ge q_*$, every $X$ as above, and every legal root position $w$,
\[
C_w\le P_*\qquad\text{or}\qquad C_w\ge q/2.
\]
The constants are independent of $q$, $w$, and the coefficients of the $f_i$.
\end{proposition}

\begin{proof}
Fix $w$ and a permutation $\pi\in\Pi(w)$. An embedding is determined by the images of the $h$ interior vertices, a point of $\A^{bh}$. The edge conditions are $a|E(U)|$ polynomial equations of degree at most $D$. Since $w$ is injective and vertices in different host classes are distinct, injectivity need only be imposed on pairs $\{x,y\}$ with $\pi(\kappa_U(x))=\pi(\kappa_U(y))$ and $\{x,y\}\not\subseteq R_U$. Let $Q_\pi$ be the set of these pairs. For each pair in $Q_\pi$, the condition $x\ne y$ is the disjoint union over $j=1,\dots,b$ of
\[
 x_1=y_1,\ \dots,\ x_{j-1}=y_{j-1},\quad x_j\ne y_j.
\]
Since $|Q_\pi|\le\binom h2+h\varrho$, these conditions partition the embeddings into at most $L=L(U,b)$ basic constructible sets; explicitly, $L\le b^{\binom h2+h\varrho}$. In a fixed piece, a pair contributes $j-1$ linear equations and one nonvanishing linear form, so the factor $b+1$ in the formula below is a uniform upper bound for the number of equations and nonvanishing conditions contributed by one pair. Thus each piece has complexity at most $(bh,m,D)$, where
\[
m=a|E(U)|+(b+1)\Bigl(\binom h2+h\varrho\Bigr).
\]
By \cref{thm:uniform-dichotomy}, each piece has at most $P$ points or at least $q/2$ points. If one piece has at least $q/2$ points, then $C_w\ge q/2$. Otherwise every piece has at most $P$ points, and summing over the at most $k!$ permutations in $\Pi(w)\subseteq S_k$ and $L$ pieces gives $C_w\le k!LP$. This proves the claim with $P_*=k!LP$.
\end{proof}

\begin{remark}\label{rem:root-fibre-constants}
Explicitly, one may take $N=bh$, $m$ as in the proof, $L=b^{\binom h2+h\varrho}$, $P_*=k!LP$, and $q_*=q_0(N,m,D)$ from \cref{thm:uniform-dichotomy}. In particular $P_*$ and $q_*$ depend only on $U,a,b,D$, as claimed.
\end{remark}

\begin{remark}
The three results enter at different points. Heintz's theorem bounds the number and degrees of the initial components, so that the constants are independent of the coefficients. B\'ezout's inequality controls the Frobenius descent. Lang--Weil, in the Cafure--Matera form, supplies the $q/2$ alternative at the terminal nodes. In Section~\ref{sec:proof}, the dichotomy is applied to the root fibres $C_w$; combined with the moment estimate, it shows that $C_w\ge q/2$ is unlikely.
\end{remark}

\section{\texorpdfstring{Proof of \cref{thm:mainresult-1}}{Proof of Theorem 1.2}}\label{sec:proof}

\subsection{The lower bound}

Let $r=a/b\in(0,1]$. By \cref{prop:rooted}, there is a finite connected bipartite rooted graph $(T,R)$ with $V(T)\setminus R\ne\emptyset$ such that $b|S|\le a e_T(S)$ for every non-empty $S\subseteq V(T)\setminus R$, and $\ex(n,T^{(p)})=O(n^{2-a/b})$ for every fixed $p$.

Let $W=\{w_1,\dots,w_{k-2}\}$ be disjoint from $V(T)$, and let $U=\mathcal S^kT$. Thus
\[
V(U)=V(T)\sqcup W,\qquad E(U)=\{e\cup W:e\in E(T)\}.
\]
Set $R_U=R\cup W$ and $A_U=V(T)\setminus R$. Put $\varrho=|R_U|$ and $h=|A_U|$. For $S\subseteq A_U$, let $e_U(S)$ be the number of edges of $U$ meeting $S$, and put
\[
\rho(U)=\frac{e_U(A_U)}{|A_U|}.
\]

Since $T$ is connected and $R$ and $V(T)\setminus R$ are both non-empty, $T$ has at least one edge. Fix its bipartition $V(T)=B_1\sqcup B_2$. The parts of $U$ are
\[
B_1,\quad B_2,\quad \{w_1\},\quad\dots,\quad \{w_{k-2}\},
\]
and we write $\kappa_U\colon V(U)\to[k]$ for the resulting part map. When $k=2$, the last $k-2$ parts are absent and $U=T$.

For $p\ge1$, the $p$-th power $U^{(p)}$ of the rooted $k$-graph $U$ is defined as for graphs: $p$ vertex-disjoint copies of $U$ are glued along the common root set $R_U$, so that they are disjoint outside $R_U$. We shall use the canonical identification
\begin{equation}\label{eq:susp-power}
	U^{(p)}=(\mathcal S^kT)^{(p)}\cong \mathcal S^k(T^{(p)}).
\end{equation}
Indeed, both $k$-graphs consist of $p$ vertex-disjoint copies of $T$, with the copies sharing the root set $R$, and with the same $(k-2)$-set $W$ added to every edge; the only difference is the order in which the disjoint union is taken. Thus they have the same vertex set up to the obvious reassociation of the disjoint union, and the same edge set. We use this identification throughout, without further comment.

\begin{lemma}\label{lem:U-parts}
Let $X$ be a $k$-partite $k$-graph with part map $\kappa_X$, and let $\varphi\colon V(U)\to V(X)$ be an embedding of $U$. Then there is a permutation $\pi\in S_k$ such that
\[
	\kappa_X(\varphi(v))=\pi(\kappa_U(v))\qquad\text{for every }v\in V(U).
\]
\end{lemma}

\begin{proof}
Choose an edge $xy\in E(T)$. The set $\varphi(\{x,y\}\cup W)$ is an edge of $X$, so its vertices lie in distinct parts; in particular the $k-2$ vertices of $\varphi(W)$ occupy distinct parts. Let $V_i,V_j$ be the two remaining parts. Every vertex $v\in V(T)$ has a neighbour $u\in V(T)$, and $\varphi(\{u,v\}\cup W)$ is an edge of $X$. Hence $\varphi(v)\in V_i\cup V_j$, and adjacent vertices of $T$ are mapped to different parts. Thus $v\mapsto\kappa_X(\varphi(v))$ is a proper $2$-colouring of the connected graph $T$ with colours $\{i,j\}$. A connected graph has a unique proper $2$-colouring up to swapping the colours, so $\varphi(B_1)$ and $\varphi(B_2)$ are contained in $V_i$ and $V_j$, respectively, in one order or the other. Together with the parts occupied by $\varphi(W)$, this determines a permutation $\pi$ with the required property. For $k=2$ the set $W$ is empty and the same argument applies.
\end{proof}

\begin{corollary}\label{cor:power-root}
Let $f\colon V(U^{(p)})\to V(X)$ be an embedding. Then $w=f|_{R_U}$ is a legal root position, and there are pairwise distinct embeddings $\varphi_1,\dots,\varphi_p\colon V(U)\to V(X)$ extending $w$, all with the same associated permutation $\pi\in S_k$.
\end{corollary}

\begin{proof}
Let $\iota_j\colon V(U)\to V(U^{(p)})$ be the inclusion of the $j$-th copy and put $\varphi_j=f\circ\iota_j$. By \cref{lem:U-parts}, each $\varphi_j$ has an associated permutation $\pi_j$. The $p$ copies share $R_U$, so the maps $\varphi_j$ agree on $R_U$; hence $w$ is a legal root position. Moreover $\pi_1=\cdots=\pi_p$: the common images of the $k-2$ vertices of $W$ determine $\pi_j$ on $k-2$ parts, and the common image of any root $r\in R$ determines it on one of the two remaining parts, hence on the last part as well. Finally, $A_U\ne\emptyset$, and the images $\varphi_j(A_U)$ are pairwise disjoint because the interiors of the $p$ copies are disjoint and $f$ is injective. Thus the $\varphi_j$ are pairwise distinct.
\end{proof}

\begin{lemma}\label{lem:U-balanced}
$U$ is balanced and $\rho(U)=b/a$.
\end{lemma}

\begin{proof}
For $S\subseteq A_U$, an edge $e\cup W$ meets $S$ if and only if $e$ meets $S$; hence $e_U(S)=e_T(S)$. Therefore $e_U(S)/|S|=e_T(S)/|S|\ge b/a$. Moreover,
\[
\rho(U)=\frac{e_U(A_U)}{|A_U|}=\frac{e_T(V(T)\setminus R)}{|V(T)\setminus R|}=\frac ba. \qedhere
\]
\end{proof}

The construction is the random algebraic method of Bukh \cite{Bukh2015}. Fix a prime power $q$. Let $V_1,\dots,V_k$ be pairwise disjoint copies of $\F_q^b$. Put
\[
s=b\varrho+a+1
\]
and choose an integer $D\ge s\,e(U)$. Choose $f_1,\dots,f_a$ independently and uniformly from the set $\mathcal P$ of polynomials of total degree at most $D$ in the $bk$ coordinates. Let $X$ be the $k$-partite $k$-graph on $V_1\cup\cdots\cup V_k$ in which a $k$-tuple $z=(z_1,\dots,z_k)\in V_1\times\cdots\times V_k$ is an edge if and only if
\[
f_1(z)=\cdots=f_a(z)=0.
\]

\begin{lemma}\label{lem:interpolation}
Let $m\le D+1$ and $q>\binom m2$. For distinct points $z_1,\dots,z_m\in\F_q^{bk}$, a uniformly random $f\in\mathcal P$ satisfies $f(z_1)=\cdots=f(z_m)=0$ with probability $q^{-m}$.
\end{lemma}

\begin{proof}
For each pair \(i\ne j\), the linear functionals \(\ell\) with \(\ell(z_i)=\ell(z_j)\) form a proper hyperplane of the dual space \((\F_q^{bk})^*\). There are \(\binom m2<q\) such hyperplanes, and each has \(q^{bk-1}\) elements. Their union has at most \(\binom m2q^{bk-1}<q^{bk}\) elements, so some \(\ell\) is injective on \(\{z_1,\dots,z_m\}\).

After an invertible linear change of coordinates, \(\ell\) is the first coordinate. Write \(t_i\) for the first coordinate of \(z_i\); the \(t_i\) are distinct. The polynomials
\[
h(t)=\sum_{j=0}^{m-1}c_jt^j
\]
form an \(m\)-dimensional subspace of \(\mathcal P\). The evaluation map \(h\mapsto(h(t_1),\dots,h(t_m))\) has the Vandermonde matrix \((t_i^j)_{1\le i\le m,\,0\le j\le m-1}\), which is invertible. Hence the evaluation map
\[
\mathcal P\longrightarrow\F_q^m,\qquad f\longmapsto(f(z_1),\dots,f(z_m)),
\]
is surjective. Therefore the values of a uniform random \(f\in\mathcal P\) at \(z_1,\dots,z_m\) are uniform in \(\F_q^m\), and the probability that all are zero is \(q^{-m}\).
\end{proof}

\begin{lemma}\label{lem:moment}
Let $s=b\varrho+a+1$. For every legal root position $w$ and every prime power $q>\binom{D+1}{2}$,
\[
\mathbb E[C_w^s]\le K,
\]
where $K=K(U,a,b,D)$ is independent of $q$, $w$, and the coefficients $f_i$.
\end{lemma}

\begin{proof}
Fix a legal root position \(w\). For a permutation \(\pi\in\Pi(w)\), let \(\mathcal C_w(\pi)\) be the set of embeddings \(\varphi:V(U)\to V(X)\) that extend \(w\) and map the parts of \(U\) according to \(\pi\), and put
\[
\mathcal C_w=\bigcup_{\pi\in\Pi(w)}\mathcal C_w(\pi),\qquad C_w=|\mathcal C_w|.
\]
We count the ordered \(s\)-tuples in \(\mathcal C_w^s\).

Before choosing the polynomials, let $\mathcal M_w(\pi)$ be the set of injective maps $V(U)\to V_1\sqcup\cdots\sqcup V_k$ extending $w$ and respecting $\pi$, without imposing edge conditions. Put $\mathcal M_w=\bigcup_{\pi\in\Pi(w)}\mathcal M_w(\pi)$; thus $\mathcal C_w\subseteq\mathcal M_w$. For a tuple $\mathbf c=(\varphi_1,\dots,\varphi_s)\in\mathcal M_w^s$, let \(\sim\) be the equivalence relation on \(A_U\times\{1,\dots,s\}\) defined by \((x,i)\sim(y,j)\) if and only if \(\varphi_i(x)=\varphi_j(y)\). Define a rooted \(k\)-graph \(H=H(\mathbf c)\) by
\[
V(H)=R_U\sqcup\bigl((A_U\times\{1,\dots,s\})/\!\sim\bigr),
\]
and
\[
E(H)=\bigl\{\psi_i(e):e\in E(U),\ 1\le i\le s\bigr\},
\]
where \(\psi_i(r)=r\) for \(r\in R_U\), \(\psi_i(x)=\overline{(x,i)}\) for \(x\in A_U\), and \(\psi_i(e)=\{\psi_i(v):v\in e\}\). Since each \(\varphi_i\) is injective and extends \(w\), the map \(\psi_i\) is injective on \(V(U)\) and fixes \(R_U\). Moreover, \(H=\bigcup_i\psi_i(U)\) and \(E(H)=\bigcup_i\psi_i(E(U))\).

The \emph{type} of \(\mathbf c\) is the tuple
\[
\tau=(H;\psi_1,\dots,\psi_s;\pi_1,\dots,\pi_s),
\]
where \(\pi_i\) is the part permutation used by \(\varphi_i\). Let \(\mathcal T\) be the set of types arising from $\mathcal M_w^s$. This set is independent of the polynomials. If $B_{sh}$ denotes the number of partitions of an $sh$-element set, then $|\mathcal T|\le B_{sh}(k!)^s$. Put $K=B_{sh}(k!)^s$.

For \(\tau=(H_\tau;\psi_1,\dots,\psi_s;\pi_1,\dots,\pi_s)\in\mathcal T\), put
\[
h_\tau=|V(H_\tau)\setminus R_U|,\qquad e_\tau=|E(H_\tau)|.
\]
Let \(N_\tau(w)\) be the number of injective maps
\[
\varphi:V(H_\tau)\longrightarrow V(X)
\]
that extend \(w\) (that is, \(\varphi(r)=w(r)\) for every \(r\in R_U\)), respect the part assignment \(\kappa_{H_\tau}\colon V(H_\tau)\to[k]\) determined by \(w\) and \(\pi_1,\dots,\pi_s\), and satisfy \(\varphi(e)\in E(X)\) for every \(e\in E(H_\tau)\). Explicitly,
\[
\kappa_{H_\tau}(r)=\kappa_X(w(r))=\pi_i(\kappa_U(r))
\quad(r\in R_U,\ 1\le i\le s),\qquad
\kappa_{H_\tau}(\overline{(x,i)})=\pi_i(\kappa_U(x)).
\]
These assignments are well-defined: the maps agree on the roots, and identified interior vertices have the same image in $X$. Injectivity is required on all of \(V(H_\tau)\), not only on the interior \(V(H_\tau)\setminus R_U\): otherwise a root image could coincide with a block image, and the maps \(\varphi\circ\psi_i\) would fail to be injective on \(V(U)\). Such a map, together with \(w\), gives the tuple
\[
(\varphi\circ\psi_1,\dots,\varphi\circ\psi_s)\in\mathcal C_w^s.
\]
Conversely, every tuple \((\varphi_1,\dots,\varphi_s)\in\mathcal C_w^s\) arises in this way from its type: the map \(\varphi\) defined by \(\varphi(r)=w(r)\) and \(\varphi(\overline{(x,i)})=\varphi_i(x)\) is well defined and injective, since \(\varphi_i(x)=\varphi_j(y)\) exactly when \((x,i)\sim(y,j)\), and it maps edges to edges. Hence
\begin{equation}\label{eq:type-decomposition}
C_w^s=\sum_{\tau\in\mathcal T}N_\tau(w).
\end{equation}

Fix \(\tau\in\mathcal T\). There are at most \(q^{b h_\tau}\) injective maps \(\varphi:V(H_\tau)\to V(X)\) extending \(w\) and respecting the part assignment (the roots are fixed by \(w\), so only the \(h_\tau\) interior images are free). Since \(D\ge s\,e(U)\), we have \(e_\tau\le s\,e(U)\le D\). For each such \(\varphi\), the \(e_\tau\) edges of \(H_\tau\) give distinct \(k\)-tuples in \(V_1\times\cdots\times V_k\): if \(\varphi(e)=\varphi(e')\) for distinct \(e,e'\in E(H_\tau)\), then injectivity of \(\varphi\) forces \(e=e'\). Applying \cref{lem:interpolation} with \(m=e_\tau\) to each of the \(a\) polynomials, we obtain
\begin{equation}\label{eq:type-expectation}
\mathbb E[N_\tau(w)]\le q^{b h_\tau-a e_\tau}.
\end{equation}

We claim that
\begin{equation}\label{eq:type-balance}
a e_\tau\ge b h_\tau.
\end{equation}
For \(1\le i\le s\), put
\[
S_i=\bigl\{x\in A_U:\psi_i(x)\notin\psi_1(A_U)\cup\cdots\cup\psi_{i-1}(A_U)\bigr\}.
\]
The sets \(\psi_i(S_i)\) are pairwise disjoint and their union is \(V(H_\tau)\setminus R_U\); hence \(\sum_i|S_i|=h_\tau\). Let
\[
E_i=\bigl\{\psi_i(e):e\in E(U),\ e\cap S_i\ne\emptyset\bigr\}.
\]
Every edge in \(E_i\) meets \(\psi_i(S_i)\), which is disjoint from \(\bigcup_{j<i}\psi_j(V(U))\); hence \(E_i\cap\psi_j(E(U))=\emptyset\) for every \(j<i\). Thus the sets \(E_i\) are pairwise disjoint. Since \(\psi_i\) is injective, \(|E_i|=e_U(S_i)\). By \cref{lem:U-balanced},
\[
|E_i|=e_U(S_i)\ge\frac ba|S_i|.
\]
Therefore
\[
e_\tau\ge\sum_{i=1}^s|E_i|
\ge\frac ba\sum_{i=1}^s|S_i|
=\frac ba h_\tau,
\]
which proves \eqref{eq:type-balance}.

Combining \eqref{eq:type-expectation} and \eqref{eq:type-balance} yields \(\mathbb E[N_\tau(w)]\le1\). Summing over \(\tau\in\mathcal T\) and using \eqref{eq:type-decomposition},
\[
\mathbb E[C_w^s]=\sum_{\tau\in\mathcal T}\mathbb E[N_\tau(w)]
\le |\mathcal T|\le K.
\]
This proves the lemma.
\end{proof}

\begin{remark}\label{rem:moment-uniform}
The type decomposition is the only place where the part permutations enter the moment estimate. Once $w$ and the permutations $\pi_1,\dots,\pi_s$ are fixed, the part of every block of $H_\tau$ is determined, so the number of choices for the interior images is at most $q^{b h_\tau}$, independently of the coefficients. The threshold $s=b\varrho+a+1$ enters the uniform bound $K=B_{sh}(k!)^s$; by \eqref{eq:type-balance}, $b h_\tau-a e_\tau\le0$, so each summand is at most $1$. Thus $K$ depends only on $U,a,b,D$, not on $q$, $w$, or the $f_i$.
\end{remark}

Let $P_*,q_*$ be as in \cref{prop:root-fibre}, and set $p_0=P_*+1$ and $F_{k,r}=U^{(p_0)}$. Take $q\ge q_*$ and $q$ large enough for the estimates below. By \cref{lem:moment} and Markov's inequality,
\[
\Pr[C_w\ge q/2]\le 2^sKq^{-s}.
\]
Since $C_w\le P_*$ or $C_w\ge q/2$, we have
\[
\Pr[C_w>P_*]\le 2^sKq^{-s}.
\]
A legal root position is determined by a permutation $\pi\in S_k$ and an injective choice of an image in the prescribed host part for each root; hence there are at most \(k!q^{b\varrho}\) legal root positions. Let \(Z\) be the number of legal root positions with \(C_w>P_*\). Then
\begin{equation}\label{eq:bad-roots}
\mathbb E[Z]\le k!2^sKq^{b\varrho-s}
=k!2^sKq^{-a-1}.
\end{equation}

Let \(M=e(X)\). By \cref{lem:interpolation} with \(m=1\), each \(k\)-tuple is an edge with probability \(q^{-a}\), so \(\mathbb E[M]=q^{bk-a}\). Let \(\mathcal B=\{Z>0\}\). Since \(Z\) is a nonnegative integer, \(\Pr[\mathcal B]\le\mathbb E[Z]\). By \eqref{eq:bad-roots}, \(\Pr[\mathcal B]\le C_0q^{-a-1}\) with \(C_0=k!2^sK\). Since \(0\le M\le q^{bk}\),
\[
\mathbb E[M\mathbf 1_{\mathcal B}]
\le q^{bk}\Pr[\mathcal B]\le C_0q^{bk-a-1}.
\]
For $q\ge2C_0$,
\begin{equation}\label{eq:many-edges}
\mathbb E[M\mathbf 1_{\mathcal B^c}]
\ge q^{bk-a}-C_0q^{bk-a-1}
\ge\frac12q^{bk-a}.
\end{equation}
The space of coefficients is finite, so the expectation in \eqref{eq:many-edges} is an average over finitely many coefficient choices; hence some choice of \(f_1,\dots,f_a\) satisfies
\begin{equation*}
Z=0,\qquad e(X)=M\ge\frac12q^{bk-a}.
\end{equation*}
For this choice, $C_w\le P_*$ for every legal root position $w$. Suppose that $X$ contained a copy of $F_{k,r}=U^{(p_0)}$, given by an embedding $f\colon V(U^{(p_0)})\to V(X)$. Let $w=f|_{R_U}$. By \cref{cor:power-root}, $w$ is a legal root position and the $p_0$ branches of the copy give $p_0$ pairwise distinct embeddings $\varphi_1,\dots,\varphi_{p_0}\colon V(U)\to V(X)$ extending $w$, all with the same associated permutation. Hence all of them are counted by $C_w$, so $C_w\ge p_0>P_*$, a contradiction. Hence $X$ is $F_{k,r}$-free.

The \(k\)-graph \(X\) has \(kq^b\) vertices and at least \(\frac12q^{bk-a}\) edges. Since \(T\) is connected and \(V(T)\setminus R\ne\emptyset\), every vertex of \(T\) is incident with an edge, and every vertex of \(W\) lies in every edge of \(U\); hence \(U\), and therefore \(F_{k,r}=U^{(p_0)}\), has no isolated vertices. For \(n\ge kq^b\), add \(n-kq^b\) isolated vertices; the resulting \(k\)-graph is still \(F_{k,r}\)-free. To cover all sufficiently large $n$, take the largest power of two $q$ with $kq^b\le n$. Then $q^b>n/(2^bk)$, and therefore
\begin{equation}\label{eq:lower-bound}
\ex_k(n,F_{k,r})
\ge\frac12q^{bk-a}
\ge\frac12(2^bk)^{-(k-a/b)}n^{k-a/b}.
\end{equation}
Thus $\ex_k(n,F_{k,r})=\Omega(n^{k-a/b})$.

\subsection{The upper bound}

By \cref{lem:suspension-upper} with $F=T^{(p_0)}$, and by the identification \eqref{eq:susp-power},
\begin{align*}
\ex_k(n,F_{k,r})
&=\ex_k(n,U^{(p_0)})
=\ex_k\bigl(n,\mathcal S^k(T^{(p_0)})\bigr)\\
&\le\binom{n}{k-2}\frac{\ex(n-k+2,T^{(p_0)})}{\binom{k}{2}}
=O\bigl(n^{k-2}\cdot n^{2-a/b}\bigr)
=O\bigl(n^{k-a/b}\bigr),
\end{align*}
where the estimate $\ex(n,T^{(p_0)})=O(n^{2-a/b})$ is from \cref{prop:rooted}. Together with \eqref{eq:lower-bound}, this proves \cref{thm:mainresult-1} for $0<r\le1$.

\subsection{The case \texorpdfstring{$r=0$}{r=0}}

Let $F_{k,0}=K_{k+1}^{(k)}$, the complete $k$-uniform hypergraph on $k+1$ vertices; its edges are all $k$-subsets of a $(k+1)$-set. The upper bound
\[
\ex_k(n,F_{k,0})\le\binom nk=O(n^k)
\]
is trivial, since every $k$-graph on $n$ vertices has at most $\binom nk$ edges.

For the lower bound, let $n\ge k$ and choose pairwise disjoint sets $P_1,\dots,P_k\subseteq[n]$ with $|P_i|=\lfloor n/k\rfloor$ for every $i$; this is possible because $k\lfloor n/k\rfloor\le n$. Let $H$ be the complete $k$-partite $k$-graph with parts $P_1,\dots,P_k$, that is, the family of all $k$-sets that meet every $P_i$ in exactly one vertex. Then
\[
e(H)=\prod_{i=1}^k|P_i|=\lfloor n/k\rfloor^k\ge\left(\frac{n}{2k}\right)^k
\qquad(n\ge 2k).
\]
Let $S\subseteq[n]$ with $|S|=k+1$. If $S$ contains a vertex outside $\bigcup_{i=1}^k P_i$, then any $k$-subset of $S$ containing that vertex is not an edge of $H$. Otherwise, two distinct vertices $x,y\in S$ lie in the same part. Since $k\ge2$, there exists $z\in S\setminus\{x,y\}$, and $S\setminus\{z\}$ is a $k$-set containing $x,y$, so it is not an edge of $H$. Thus $H$ contains no copy of $K_{k+1}^{(k)}$. Therefore
\[
\ex_k(n,F_{k,0})\ge\lfloor n/k\rfloor^k=\Omega_k(n^k),
\]
and hence $\ex_k(n,F_{k,0})=\Theta(n^k)$.

\section{A two-edge obstruction}\label{sec:obstruction}

For $0\le t\le k-1$, let $Q_{k,t}$ be the $k$-graph consisting of two $k$-edges that intersect in exactly $t$ vertices.

\begin{lemma}\label{lem:twoedge}
For fixed $k,t$, $\ex_k(n,Q_{k,t})=\Omega_k(n^{\max\{t,k-t-1\}})$.
\end{lemma}

\begin{proof}
Fix a $(t+1)$-set $C$ and take all $k$-sets containing $C$. Any two edges meet in at least $t+1$ vertices, so the family is $Q_{k,t}$-free and has $\binom{n-t-1}{k-t-1}=\Omega_k(n^{k-t-1})$ edges.

For $t\ge1$, take a maximal family $\mathcal M\subseteq\binom{[n]}k$ whose distinct members meet in at most $t-1$ vertices. It is $Q_{k,t}$-free. By maximality, every $k$-set $A$ satisfies $|A\cap Q|\ge t$ for some $Q\in\mathcal M$. Hence every $k$-set contains a $t$-set that is contained in some member of $\mathcal M$. There are at most $|\mathcal M|\binom kt$ such $t$-sets, and each is contained in $\binom{n-t}{k-t}$ $k$-sets. Thus
\begin{equation*}
	\binom nk\le|\mathcal M|\binom kt\binom{n-t}{k-t},
	\qquad
	|\mathcal M|\ge\frac{\binom nt}{\binom kt^2}=\Omega_k(n^t).
\end{equation*}
For $t=0$, a single edge gives a constant lower bound. Taking the larger of the two constructions proves the lemma.
\end{proof}

\begin{proof}[Proof of \cref{prop:obstruction}]
If $F$ has at least two edges, choose two edges meeting in $t$ vertices. Every $Q_{k,t}$-free $k$-graph is $F$-free, so, by \cref{lem:twoedge},
\begin{equation*}
	\ex_k(n,F)\ge\ex_k(n,Q_{k,t})=\Omega_k(n^{\max\{t,k-t-1\}})=\Omega_k(n^{\lfloor k/2\rfloor}).
\end{equation*}
Since $r>\lceil k/2\rceil$, we have $k-r<\lfloor k/2\rfloor$, so $\ex_k(n,F)=\Omega(n^{\lfloor k/2\rfloor})$ contradicts $\ex_k(n,F)=\Theta(n^{k-r})$.

If $F$ has no edges, then every $k$-graph on at least $|V(F)|$ vertices contains a copy of $F$, so $\ex_k(n,F)=0$. If $F$ has exactly one edge, then for $n\ge|V(F)|$ an $F$-free $k$-graph has no edges, so again $\ex_k(n,F)=0$. In either case $\ex_k(n,F)=O(1)$, which is not $\Theta(n^{k-r})$ since $k-r\ge1$.
\end{proof}

\section{Concluding remark}\label{sec:concluding-remark}

\Cref{thm:mainresult-1} settles the endpoint $r=1$ of the interval $[1,k/2]$. For $1<r\le k/2$ we have $k-r\ge\lfloor k/2\rfloor$, so the lower bound $\Omega(n^{\lfloor k/2\rfloor})$ of \cref{prop:obstruction} does not exclude the desired order of magnitude. We propose the following extension.

\begin{conjecture}\label{conj:single-hypergraph}
For every integer $k\ge3$ and every rational $r\in[1,k/2]$, there exists a $k$-uniform hypergraph $F$ such that
\[
\ex_k(n,F)=\Theta(n^{k-r}).
\]
\end{conjecture}

The suspension inequality used here transfers graph exponents only to $0\le r\le1$, so \cref{conj:single-hypergraph} requires a different upper-bound argument. For odd $k$, the additional interval $k/2<r\le\lceil k/2\rceil$ is neither covered by this conjecture nor excluded by \cref{prop:obstruction}.

\section*{Declaration of competing interest}
The authors declare no competing interests.

\section*{Acknowledgements}
We are grateful to Thomas Bloom for establishing and curating the Erd\H{o}s Problems website, which served as an invaluable resource and source of inspiration for this research.

\section*{Declaration of AI usage}
This work was carried out with the assistance of DeepSeek-V4.1-Flash 1017. The AI system was used for discussions of the mathematical argument, for language editing and polishing of the text, and for the Lean 4 formalization of \cref{thm:mainresult-1}.

\Cref{thm:mainresult-1} has been formalized in Lean 4\footnote{The formalized proof is available at \url{https://github.com/cnudt/erdos-problem-571-hypergraph}.}, conditional on two axioms: the literature input, consisting of \cref{prop:rooted}, the suspension inequality of \cref{lem:suspension-upper}, and the identification \eqref{eq:susp-power}; and the uniform dichotomy of \cref{thm:uniform-dichotomy}. Apart from these, the argument is formalized; in particular \cref{prop:root-fibre} is proved rather than assumed. \cref{prop:obstruction} is not part of the formalization. The authors have checked the argument and the final text and take full responsibility for the mathematical content.

\end{document}